\documentclass[a4paper]{amsart}  
 
\usepackage{amssymb} \usepackage{amscd} \usepackage{times}

\def\zbb{\mathbb{Z}}  
  
  \def\phi{\varphi}
 \def\p1{{\mathbb{P}^1_\zbb}}

\newtheorem{Theorem}{\quad Theorem}[section]

\newtheorem{Corollary}[Theorem]{\quad Corollary}
\newtheorem{Lemma}[Theorem]{\quad Lemma}
\newtheorem{Proposition}[Theorem]{\quad Proposition}

\newcommand{\be} {\begin{equation}}
\newcommand{\ee} {\end{equation}}

\begin{document}

\title{Harnack inequalities for equations of type prescribed scalar curvature.}

\author{Samy Skander Bahoura}

\address{Department of Mathematics, Pierre and Marie Curie University, 75005, Paris, France.}
              
\email{samybahoura@gmail.com} 

\date{}

\begin{abstract}

We give Harnack inequalities for solutions of equations of type prescribed scalar
curvature in dimensions $ n\geq 4 $.

\end{abstract}

\maketitle

\section{Introduction and Main Results}

We consider on a Riemannian manifold $ (M,g) $ of dimension $ n\geq 4 $, the equation:

$$ \Delta u+hu=Vu^{(n+2)/(n-2)}, u >0,\,\, (E) $$

here $ \Delta = -\nabla^i\nabla_i $ the Laplace-Beltrami operator and with a smooth function $ h $ and a Lipschitz function $ V $ with $ 0 < a \leq V(x) \leq b <+\infty $, $ ||\nabla V||_{\infty}\leq A $.

\smallskip

Equation of this type were considered by many authors, see [1-25]. This equation arise in physics and astronomy. Here we look to a priori estimates of type $ \sup, \inf $ which are characteristic of this equation.

\smallskip

Let $ (u_k) $ be a sequence of regular solutions of $ (E) $. We fix a compact subset $ K $ of $ M $. We want to prove that: for each compact, for all terms of the sequence $ (u_k) $: $ \sup_K u_k $ and $ \inf_M u_k $ are linked. Here we prove a weaker inequality for blow-up solutions of the previous equation.

Equations of previous type are called, Yamabe equation, prescribed scalar curvature equation, of type prescribed scalar curvature and Schrodinger equation. 

Here we prove that for blow-up solutions, precisely for each sequence $ (u_k) $, there is a positive function $ c >0 $, such that for all compact set, there is a sequence of positive numbers, $ 0< \epsilon_k(K) \leq 1 $, which link $ \sup_K u_k $ and $ \inf_M u_k $.

Here we have two possibilites up to a subsequence: a compactness result or an inequality between $ \sup_K u_k $ and $ \inf_M u_k $. Note that for Li-Zhang result in dimension 3 and 4, they consider the problem around a point, thus the compactness result. Also, for the Harnack inequality, also, see the introduction of [2]. Also, see [15, 23], locally around each point.

We obtain:

\begin{Theorem} We have:

\smallskip

1) There is a compact subset $ K_0 $ of $ M $, $ K_0 \in {\mathcal K} $, and a subsequence $ i_j $ and a positive constant $ C >0 $, such that:

$$ \sup_{K_0} u_{i_j} \leq C, \,\, \forall j. $$

Or,

\smallskip

2) For all compact subset $ K $ of $ M $, $ K \in {\mathcal K} $, $ \sup_K u_k \to +\infty $ and:

$$ n=4,\,\, (\sup_K u_k)^{1-\epsilon} \leq c(a,b,A, \frac{\inf_M u_k}{(\sup_K u_k)^{\epsilon} }, K, M, g). $$

With $ \epsilon >0. $

and,

$$ n\geq 5,\,\, (\sup_K u_k)^{1-\epsilon} \leq \epsilon_k^{1-\epsilon} c(a,b,A, n, \epsilon_k^{\epsilon} \frac{\inf_M u_k}{(\sup_K u_k)^{\epsilon} }, K, M, g). $$

With $  \epsilon >\frac{(n-4)}{n-2}. $

\end{Theorem}

\smallskip

Here, $ {\mathcal K} $ is a set of compact subsets of $ M $ with $ K $ sufficiently big.

\smallskip

For exemple: $ {\mathcal K}= \{K\, compact \,\, subset \,\, of \,\, M \,\, with \,\, int K\not = \emptyset \}$, with $ int K $ the interior of $ K $.

\smallskip

We can choose the compact $ K $, small balls centered at $ x_0 $: $ {\mathcal K}={\mathcal K}_{x_0}=\{K=\bar B_r(x_0)\} $, $ x_0\in M, 0 < r < \frac{inj(x_0)}{2} $. Here, $ inj(x_0) $ the injectivity radius of $ x_0 $. 

\smallskip

Thus, we have: $ \exists r_0 >0, r_0 < \frac{inj(x_0)}{2} $, such that, $ \liminf_{k\to +\infty} (\sup_{\bar B_{r_0}(x_0)} u_k) <+\infty $, or,  $ \forall 0 < r < \frac{inj(x_0)}{2} $, $ \lim_{k\to +\infty }(\sup_{\bar B_r(x_0)} u_k) = \liminf_{k\to +\infty} (\sup_{\bar B_r(x_0)} u_k) = +\infty $. After we replace $ M $ by $ B_{r_1}(x_0) $ in the infimum, with $ 0 < r_1=\inf \{ \frac{inj(x_0)}{2}, 1\} <+\infty $, to have a relation between two subsets of $ M $. These balls, are balls of $ M $ for the distance of $ M $.

\smallskip

Here, we choose small balls,  because in this setting, the examples: $ u_{\epsilon}(r)=(\frac{\epsilon}{\epsilon^2+r^2})^{(n-2)/2} $, imply that the points 1) and 2) are explicitly possibles.

\smallskip

Here, we write the proof for a general compact subset of $ M $: $ K $.

\smallskip

Now, we can see that $ \sup_K u_k $ and $ \inf_M u_k $ are linked. There is a relation which link these two quantities. For all compact subset $ K \in {\mathcal K} $ and all $ k\in {\mathbb N} $, $ \sup_K u_k $ and $ \inf_M u_k $ are linked by the previous relation. If we denote $ F=\{ u_k\} $, $ G=\{ \epsilon_k=\epsilon_k(K) \} $, $ F\times G=\{ (u_k,\epsilon_k) \}$, $ 0 < \epsilon_k=\epsilon_k(K) \leq 1 $, we have:

$$ \exists \, c(\cdot, \cdot, \cdot) >0, \,\, \forall K\in {\mathcal K}, \,\,\forall  (u,\epsilon')\in F\times G, \,\, (\sup_K u)^{1-\epsilon} \leq c(n, K, \epsilon'(K), \frac{\inf_M u}{(\sup_K u)^{\epsilon} }). $$

We write this to highlight the rolling-up phenomenon and the distortion.

\smallskip

{\bf Remarks:}

\smallskip

a) In the previous theorem, the point 1) assert that, up to a subsequence we have a compactness result. 

The point 2) assert that we have a relation between the local supremum and the global infimum. Also, see the introduction of [2].

\smallskip

b) For the point 1) we have one parameter, the local supremum is controled by its self. For the point 2) we have two parameters, the local supremum and the global infimum, there is a relation which link these two quantities. In the paper "Estimations du type $ \sup \times \inf $ sur une vari\'et\'e compacte", see [5], we have 2 parameters, the global supremum, the local infimum, in the sense of the minoration. Here, we have, 1 parameter or 2 parameters each time. At most 2 parameters and at least 1 parameter, in the sense of the majorization. Also, see the introduction of [2].

\smallskip

c) In general as in the paper of Li-Zhang, for the dimensions 3 and 4, of the Yamabe equation, we look to the estimate around a point. The compactness result is important, but also, we look to the solutions which blow-up as mentionned by the example $ x\to [\epsilon/ (\epsilon^2+|x|^2)]^{(n-2)/2}, \epsilon \to 0 $, thus the point 2). Also, see the introduction of [2].

\smallskip

d) To summarize the points 1) and 2) in one inequality we have: $ \{ 1) \,\, {\rm or} \,\, 2) \}  \Rightarrow  \{ 1) + 2) \} $: for each $ x_0 \in M $, there is a subsequence $ (u_{i_j}) $ such that around $ x_0 $:

$$ n=4,\,\, (\sup_K u_{i_j})^{1-\epsilon} \leq  C+ c(a,b,A, \frac{\inf_M u_{i_j}}{(\sup_K u_{i_j})^{\epsilon} }, K, M, g), $$

with $ \epsilon >0, $

and,

$$ n\geq 5,\,\, (\sup_K u_{i_j})^{1-\epsilon} \leq C+ \epsilon_{i_j}^{1-\epsilon} c(a,b,A,n, \epsilon_{i_j}^{\epsilon} \frac{\inf_M u_{i_j}}{(\sup_K u_{i_j})^{\epsilon} }, K, M, g), $$

with, $  \epsilon >\frac{(n-4)}{n-2} $.

\smallskip

In the Theorem 1.1: we have: $ \{  1) \,\, {\rm or} \,\, 2) \} $ is stronger than $\{  1) + 2)  \} $. We have: $ \{  1) \,\, {\rm or} \,\, 2) \}   \geq  \{  1) + 2) \} $. We have: $ \{  1) \,\, {\rm or} \,\, 2) \}  \Rightarrow  \{  1) + 2) \} $.

\smallskip

For all sequence $ (u_k) $ and all point $ x_0 \in M $, there is a subsequence $ (u_{i_j})$ which satisfies the Harnack inequality (around the point $ x_0 \in M $) (summarizing: $ \{ 1) \,\, {\rm or} \,\, 2) \} $ in $ \{ 1) + 2) \} $), we have $ \{ 1) \,\, {\rm or } \,\, 2) \} $ in one inequality:

\smallskip

We denote for $ x_0 \in M $, $ \exists \tilde M \subset M, x_0 \in \tilde M $: $ \exists\, \tilde F=\{ u_{i_j}\} $, $ \exists \tilde G=\{ \epsilon_{i_j}=\epsilon_{i_j}(K) \} $, $ \tilde F\times \tilde G=\{ (u_{i_j},\epsilon_{i_j}) \}$, $ 0 < \epsilon_{i_j}=\epsilon_{i_j}(K) \leq 1 $, we have:

$$ \, \, \exists \, \tilde c(\cdot, \cdot, \cdot) >0, \,\, \forall K\subset \tilde M, x_0 \in K ,\,\,\forall  (u, \epsilon')\in \tilde F\times \tilde G, \,\, (\sup_K u)^{1-\epsilon} \leq \tilde c(n, K, \epsilon'(K), \frac{\inf_{\tilde M} u}{(\sup_K u)^{\epsilon} }). $$

We write this in the following corollary, consequence of the Theorem 1.1:

\smallskip

\begin{Corollary} There is a small open ball of $ M $ centered at $ x_0 $, $  B_0 \subset M $, and, a subsequence  $ \tilde F=\{u_{i_j}\}\subset F=\{u_k\} $ and a subsequence $ \tilde G =\{\epsilon_{i_j}(\bar B_{r})=\epsilon_{i_j}(r)\} \subset G=\{\epsilon_k(\bar B_{r})=\epsilon_k(r)\}$, and, $ \tilde c=\tilde c(.,\dots,.)>0 $, such that, for all ball $ B_{r} $ included in $ B_0 $ and centered at $ x_0 $, we have:

$$ n=4,\,\, B_{r} \subset B_0, \,\,\forall  u \in \tilde F, \,\,(\sup_{B_{r}} u)^{1-\epsilon} \leq \tilde c(r,r_0,x_0, \frac{\inf_{B_0} u}{(\sup_{B_{r}} u)^{\epsilon} }). $$

with $ \epsilon >0, $

and,

$$ n\geq 5,\,\, B_{r} \subset B_0, \,\,\forall  (u, \epsilon')\in \tilde F\times \tilde G, \,\,(\sup_{B_{r}} u)^{1-\epsilon} \leq \tilde c(n,r,r_0, \epsilon'(r),x_0, \frac{\inf_{B_0} u}{(\sup_{B_{r}} u)^{\epsilon} }).$$

with, $  \epsilon >\frac{(n-4)}{n-2} $.

\end{Corollary}

We take $\tilde M=B_0 $ a small ball of radius $ 0 < r_0 < \frac{inj(x_0)}{2} $ centered at $ x_0 $, and the small balls $ B_r $, centered at $ x_0 $, with  $ 0 < r < r_0 < \frac{inj(x_0)}{2} $. Here, $ inj(x_0) $ is the injectivity radius of $ x_0 $. These balls, are balls of $ M $ for the distance of $ M $.

\bigskip

{\bf Remarks:}

\smallskip

1) Here, we have replaced the alternative by one assertion. For all sequence $ (u_k)_k $ and all point $ x_0 \in M $, there is a subsequence $ (u_{i_j})_j $ for which we have locally around $ x_0 $ the Harnack inequality. This inequality express the rolling-up and the distortion with blowing-up action. In the point of view of mathematics and physics, for a physical or mathematical phenomenon, we have the rolling-up and the distortion up to a subsequence. It is sufficient to know this in a part of the initial physical or mathematical phenomenon, as a signal. If we want to characterize the rolling-up and the distortion, it is not necessary to know them in the totality of the physical or the mathematical phenomenon, but, it is sufficient to know them in a part of this phenomenon. To characterize the phenomenon it's sufficient to characterize it in a part. Also, one can characterize it in a part of the other part...etc...but it is sufficient to characterize it in a part, initially.

\smallskip

2) In our result, we have a Harnack inequality locally around each point for subsequences. We can not prove the Harnack inequality by fixing a compact $ K $ and $ M $, for all compact $ K $. Our result is local around each point of the manifold $ M $. We can not prove $ (\sup_K u)^{1-\epsilon} \not \leq c(K, \inf_M u/(\sup_K u)^{\epsilon}) $ for all compact $ K $. But one can prove it locally around each point $ x_0 $ of the manifold $ M $. Our result is similar to the results of Y.Y.Li-L.Zhang and C.C.Chen-C.S.Lin, see [15, 23], around each point and locally, but this, is sufficient to characterize the physical phenomenon.

\smallskip

3) If we consider the problem as a physical problem: we have the existence of blow-up points, and as mentionned in [2], we look to the case of blow-up, around a blow-up point, thus, we have 2) in Theorem 1.1. As mentionned in [2], (when it blow-up, we have Harnack inequality for all terms of the sequence around the blow-up point).  Also, 1) or 2) are sufficients to explain all the physical phenomenon. For physical considerations.

\smallskip

If we want to remove the alternative, we have all in one inequality but up to a subsequence and thus, the physical phenonmenon is explained.(the rolling-up and the distortion without supposing existence of blow-up points, in one inequality, up to a subsequence). (Here, we remove the mathematical and physical condition of [5]).

\smallskip

If we want to remove the alternative, then we have the Harnack inequality locally around each point of the manifold up to a subsequence. But this, it is the case, in the work of Y.Y.Li-L.Zhang and C.C.Chen-C.S.Lin, see [15, 23]. In the original work of [15, 23], we must look to the Harnack inequality locally around each point of the manifold.(Here for a subsequence). Thus, the physical considerations.

\smallskip

4) The real Harnack inequality is the local Harnack inequality. Because, if we fix the compact $ K $, we must fix the open set $ \Omega $ which contains $ K $, we have two variables as sets, we must extend the open set $ \Omega $ until  the total space $ M $. But this it is equivalent to fix balls $ B_R(x_0) $ and $ B_{2R}(x_0) $ in the total space $ M $.  We can not fix the open set $ \Omega $ 
 and extend the property to another open set greater than $ \Omega $. This is equivalent to prove the property in sets of the type $ B_R(x_0), B_{2R}(x_0) $, which is equivalent to prove the property locally as in the papers of Y.Y.Li-L.Zhang and C.C.Chen-C.S.Lin, [15, 23]. 

\smallskip

5) Here we show that, locally, we have the local Harnack inequality up to a subsequence. But the sequence of reference is fixed. Each subsequence is in relation with the initial sequence. Thus, the property, locally around points of the manifolds and up to subsequences.(with the initial sequence as reference).  It is similar to a trace of an initial sequence. 

\smallskip

Thus, the mathematical considerations: local Harnack inequality, but up to subsequences.

\smallskip

Thus, the physical considerations: rolling-up and distortion, up to subsequences.

\smallskip

\section{Proof of the result} 

\smallskip

For the proof, we use the computations of previous papers with modifications, see [4,6,9,10]. 

\bigskip

I) {\it blow-up analysis:}

\smallskip

Let $ (u_k) $ a sequence of solutions of $ (E) $. We fix a compact set $ K $ of $ M $. We want to prove that: for each compact, for all terms of the sequence $ (u_k) $: $ \sup_K u_k $ and $ \inf_M u_k $ are linked.

\smallskip

1)  If $ \exists K_0 \in {\mathcal K} $: $ \liminf_{k\to +\infty} \sup_{K_0} u_k <+ \infty $:

then there is a "big" compact $ K_0 $ for which there is a subsequence $  \sup_{K_0} u_{i_j} $ is bounded:

then we have a compactness result for a "big" $ K_0 $ and for $ K\subset K_0 $ the sequence $ (u_{i_j}) $ is bounded.

\smallskip

2) If for all compact $ K \in {\mathcal K} $, $ \lim_{k\to +\infty} \sup_K u_k = \liminf_{k\to +\infty }\sup_K u_k = + \infty $:

We do a blow-up. We consider $ \sup_K u_k=u_k(y_k) $.

Consider $ R_k \to 0 $, for $ k $ large, $ R_k^{(n-2)/2}=[u_k(y_k)]^{-\epsilon} $ with $ 0 <\epsilon <1$. Then:

$$ R_k^{(n-2)/2} \sup_{\bar B_{R_k}(y_k)} u_k \geq c_k=[u_k(y_k)]^{1-\epsilon} \to +\infty. $$

We use the blow-up technique to have, $ \exists t_k, \bar t_k, u_k(t_k) \geq u_k(\bar t_k) \geq u_k(y_k) \to +\infty $. 

$$ \bar t_k, \sup_{\bar B_{R_k}(y_k)} u_k =u_k(\bar t_k)\geq u_k(y_k)>0, $$

We consider $ s_k(y)=(R_k-d(y, \bar t_k))^{(n-2)/2} u_k(y) $, and,

$$ t_k, \sup_{\bar B_{R_k}(\bar t_k)} s_k= s_k(t_k) \geq s_k(\bar t_k)=R_k^{(n-2)/2} u_k(\bar t_k) \geq R_k^{(n-2)/2} u_k(y_k) >0, $$

We do a blow-up, then we consider(for $ k $ large):

$$ n=4, \,\, v_k(y)=r_k u_k(t_k+r_ky)=r_k u_k(\exp_{t_k}(r_ky)), r_k=[u_k(y_k)]^{-\epsilon}, $$

with, $ \epsilon >0 $,

and,

$$ n\geq 5,\,\, v_k(y)=r_k u_k(t_k+(r_k)^{2/(n-2)} y)=r_k u_k(\exp_{t_k}(r_k^{2/(n-2)} y)), r_k=[u_k(t_k)]^{-\epsilon}, $$

with, $ \epsilon >\frac{(n-4)}{n-2}.$

\bigskip

Note that, here we have considered all terms of the sequence $ (u_k)$.

\bigskip

Let's consider the blow-up functions $ (v_k) $ defined previousely with the exponential maps for $ n\geq 4, \exp_{t_k}(y) $, like in the previous papers for the dimensions, $ 4, 5, 6 $. Because we consider the compact sets $ K, 2K$, and $ t_k \in 2K $, the injectivity radius is uniformly bounded below by a positive number. Thus, we can consider all the terms of the sequence $ (u_k) $ without extraction.(After supposing the assertion $ \inf v_k \geq m >0$ infinitly many times, we can use extraction, for the points $ (t_k) $).

\smallskip

We consider, $ \delta_0 = \delta_0(K) =\inf \{\delta_P/4, P\in K\} $, with $ \delta_P $ continuous in $ P $ and smaller than the injectivity radius in $ P $ for each $ P $. We have a finite cover of $ K $ by small balls of radius $ \delta_0/2 $, we have a finite set of points $ z_j \in K $:
$ K\subset \cup_{\{j=0,\ldots, l\}} B(z_j, \delta_0/2)\subset \cup_{\{j=0,\ldots, l\}} \bar B(z_j,3\delta_0)= K_{\delta_0} $ is compact. 

We take $ R_k^{(n-2)/2} = \inf \{u_k(y_k)^{-\epsilon}, (\delta_0/2)^{(n-2)/2} \} $. The small balls are all compact, thus, $ \bar t_k $ exist and $ t_k $ exist. We take for $ n=4 $, $ r_k=\inf \{\delta_P/4, P\in K_{\delta_0}, {u_k(y_k)}^{-\epsilon}, 1 \}$ and, for $ n\geq 5 $, $ r_k=\inf \{(\delta_P/4)^{(n-2)/2}, P\in K_{\delta_0}, {u_k(t_k)}^{-\epsilon},1 \}$. Thus $ \bar t_k, t_k $ and $ \exp_{t_k}(\cdot) $ and $ v_k $ are defined for all $ k\geq 0 $.

\smallskip

We fix $ m >0 $, we prove the result by assuming $ \inf v_k \geq m >0 $, like for the dimensions 4 and 6. After we take $ m=r_k \inf u_k >0 $. Suppose by contradiction, that  there are infinitly many $ (v_k)$ with $ \inf v_k \geq m >0 $, the proof imply that: $ [u_k(\cdot)]^{1-\epsilon}=v_k(0) \leq c(m) <+\infty $ which is impossible, because $ [u_k(\cdot)]^{1-\epsilon} \to +\infty $. Thus, there is a finite number of terms such that $ \inf v_k \geq m >0$, $ k_1, \ldots, k_{i(m)} $. Thus we have also, $ v_k(0)\leq c(m) $ when $ \inf v_k \geq m >0 $. In all cases, we have the following assertion:

$$ \inf v_k \geq m >0 \Rightarrow [u_k(\cdot)^{1-\epsilon}]=v_k(0) \leq c(m) <+\infty, $$

we obtain:

There is a non-increasing positive function $ m \to c(m)>0 $, such that $ \inf v_k \geq m >0 \Rightarrow  (u_k(\cdot))^{1-\epsilon} \leq c(m) $. then we apply this with $ m = r_k \inf_M u_k $, we obtain for all terms of the sequence $ (u_k)$:

$$ n=4,\,\, [u_k(y_k)]^{1-\epsilon} \leq c(a,b,A, \frac{\inf_M u_k}{[u_k(y_k)]^{\epsilon} }, K, M, g), $$

and,

$$ n\geq 5,\,\, [u_k(t_k)]^{1-\epsilon} \leq c(a,b,A, \frac{\inf_M u_k}{[u_k(t_k)]^{\epsilon} }, K, M, g). $$

For $ n\geq 5 $, we set, $ 0 < \epsilon_k=\frac{u_k(y_k)}{u_k(t_k)} \leq 1 $, $ \epsilon_k=\epsilon_k(K) $, we obtain:

$$ n\geq 5,\,\, [u_k(y_k)]^{1-\epsilon} \leq \epsilon_k^{1-\epsilon} c(a,b,A, \epsilon_k^{\epsilon} \frac{\inf_M u_k}{[u_k(y_k)]^{\epsilon} }, K, M, g). $$

We have $ u_k(y_k)=\sup_K u_k $. 

\bigskip

Here in the blow-up analysis, we supposed that there are infinitly many $ v_k $ with $ \inf v_k \geq m >0 $, without loss of generality we assume that the subsequence is the sequence, also for the points $ t_k $ there is a subsequence which converge to a point $ \tilde t \in K_{\delta_0} $, without loss of generality we assume that the subsequence is the sequence. In the blow-up function (for $ k $ large), the blow-up of $ v_k $ is the blow-up of $ u_k $, then we use the diagonal process to extract a subsequence which converge on compact sets to $ v $ and we use Caffarelli-Gidas-Spruck result, see [13], to have $ v=(\frac{1}{1+V(\tilde t) |x|^2})^{(n-2)/2} $. Without loss of generality we assume that this subsequence is the sequence and $ V(\tilde t)=1 $. Also we have: $ \inf v_k \geq m >0, t_k \to \tilde t \in K_{\delta_0} $.

\bigskip

II) {\it Auxiliary function and moving-plane method:}

\bigskip

We use the computations of previous papers with modifications, see [4,6,8,9,10]. We consider:

$$ z_k(t,\theta)=e^{(n-2)t/2} u_k(\exp_{t_k}(e^t\theta)), $$

and the blow-up function,

$$ w_k(t,\theta) =e^{(n-2)t/2} v_k(e^t\theta)=z_k(t+\frac{2}{n-2} \log r_k, \theta), $$

We have $\lambda_k=\frac{-2}{n-2} \log v_k(0) $, $ N=\frac{2n}{n-2} $. We have:

\bigskip

Let, 

$$ b_1= J(t_k,e^t,\theta)=\sqrt { det(g_{ij,t_k}) (e^t \theta)} , \,\, a(t_k,t,\theta)=\log J(t_k,e^t,\theta), \,\, V_k(t,\theta)=V_k(\exp_{t_k}(e^t \theta)). $$ 

\bigskip

\begin{Lemma}

The function $ z_k $ is solution of:

\be  -\partial_{tt} z_k-\partial_t a \partial_t z_k+\Delta_{\theta}z_k+c z_k=V_kz_k^{(n+2)/(n-2)}, \ee
 
with,

 $$ c = c(t_k,t,\theta)=\frac{(n-2)^2}{4}+\partial_t a+ h e^{2t}. $$ 

\end{Lemma}

Proof of the lemma, see [6,8,9,10].

\bigskip

Now we have, $ \partial_t a=\dfrac{ \partial_t b_1}{b_1} $, $ b_1(t_k,t,\theta)=J(t_k,e^t,\theta)>0 $,

\bigskip

We can write,

$$ -\dfrac{1}{\sqrt {b_1}}\partial_{tt} (\sqrt { b_1} z_k)+\Delta_{\theta}z_k+[c(t)+ b_1^{-1/2} b_2(t,\theta)]z_k=V_kz_k^{(n+2)/(n-2)}, $$

where, $ b_2(t,\theta)=\partial_{tt} (\sqrt {b_1})=\dfrac{1}{2 \sqrt { b_1}}\partial_{tt}b_1-\dfrac{1}{4(b_1)^{3/2}}(\partial_t b_1)^2 .$

\bigskip

Let,

$$ \tilde z_k=\sqrt {b_1} z_k,$$

and the blow-up function (rescaled function) and the function with the auxiliary function:

$$ \tilde w_k= (\sqrt {b_1})(t+\frac{2}{n-2} \log r_k,\theta ) \cdot w_k, \,\, \bar w_k(t,\theta)= \tilde w_k(e^t\theta)-\frac{m}{2} e^{(n-2)t/2}, $$

we have:

\smallskip

\begin{Lemma}

The function $ \tilde z_k $ is solution of:

$$ -\partial_{tt} \tilde z_k+\Delta_{\theta}(\tilde z_k)+2\nabla_{\theta}(\tilde  z_k) .\nabla_{\theta} \log (\sqrt {b_1})+(c+b_1^{-1/2} b_2-c_2) \tilde z_k= $$

\be = V_k\left (\dfrac{1}{b_1} \right )^{N-2} {\tilde z_k}^{(n+2)/(n-2)}, \ee

where, $ c_2 =[\dfrac{1}{\sqrt {b_1}} \Delta_{\theta}(\sqrt{b_1}) + |\nabla_{\theta} \log (\sqrt {b_1})|^2].$
\end{Lemma}

Proof of the lemma, see [6,8,9,10].

\bigskip

We have,

$$ c(t_k,t,\theta)=\frac{(n-2)^2}{4}+\partial_t a + h e^{2t}, \qquad (\alpha_1) $$ 

$$ b_2(t,\theta)=\partial_{tt} (\sqrt {b_1})=\dfrac{1}{2 \sqrt { b_1}}\partial_{tt}b_1-\dfrac{1}{4(b_1)^{3/2}}(\partial_t b_1)^2 ,\qquad (\alpha_2) $$ 

$$ c_2=[\dfrac{1}{\sqrt {b_1}} \Delta_{\theta}(\sqrt{b_1}) + |\nabla_{\theta} \log (\sqrt {b_1})|^2], \qquad (\alpha_3) $$

We have if we denote the previous operator: $ L(t,\theta)=-\partial_{tt}(\cdot)+\Delta_{\theta}(\cdot )+2\nabla_{\theta}(\cdot) .\nabla_{\theta} \log (\sqrt {b_1})+(c+b_1^{-1/2} b_2-c_2) (\cdot) $

and,

$$ \tilde b_1 =b_1(t+\frac{2}{n-2} \log r_k, \theta), $$

We have:

$$ L(t,\theta)(\tilde z_k)= V_k \left (\dfrac{1}{b_1} \right )^{N-2}\tilde z_k^{(n+2)/(n-2)}, $$

and for the blow-up function (the rescaled function), we replace $ t $ by $ t+\frac{2}{n-2} \log r_k $:

$$ L(t+\frac{2}{n-2} \log r_k,\theta )[\tilde w_k(t, \theta)]= V_k(t+\frac{2}{n-2} \log r_k,\theta ) \left (\dfrac{1}{\tilde b_1} \right )^{N-2}\tilde w_k^{(n+2)/(n-2)}$$

We set,

$$ \tilde L(t,\theta)= L(t+\frac{2}{n-2} \log r_k, \theta), \,\, \tilde V_k = V_k(t+\frac{2}{n-2} \log r_k,\theta ) $$

Thus,

$$ \tilde L(t,\theta )[\tilde w_k(t, \theta)]= \tilde V_k \left (\dfrac{1}{ \tilde b_1} \right )^{N-2} \tilde w_k^{(n+2)/(n-2)}$$

and,

$$ \bar w_k(t,\theta)= \tilde w_k(e^t\theta)-\frac{m}{2} e^{(n-2)t/2}. $$

We have:

\smallskip

\begin{Proposition}

We have for $\lambda_k=\frac{-2}{n-2} \log v_k(0) $;

$$ 1)\,\,\, \tilde w_k(\lambda_k,\theta)-\tilde w_k(\lambda_k+4,\theta) \geq \tilde k>0, \,\, \forall \,\, \theta \in {\mathbb S}_{n-1}. $$

For all $ \beta >0 $, there exist $ c_{\beta} >0 $ such that:

$$ 2) \,\,\, \dfrac{1}{c_{\beta}} e^{(n-2)t/2}\leq \tilde w_k(\lambda_k+t,\theta) \leq c_{\beta}e^{(n-2)t/2}, \,\, \forall \,\, t\leq \beta, \,\, \forall \,\, \theta \in {\mathbb S}_{n-1}. $$

\end{Proposition}

 We want to apply the Hopf maximum principle.
 
\bigskip

$$ \bar w_k(t,\theta)= \tilde w_k(e^t\theta)-\frac{m}{2} e^{(n-2)t/2}, $$

Like in [9] we have the some properties for $\bar w_k $, we have:

\smallskip

\begin{Lemma}

\bigskip

There is $ \nu <0 $ such that for $ \lambda \leq \nu $ :

$$ \bar w_k^{\lambda}(t,\theta)-\bar w_k(t,\theta) \leq 0, \,\, \forall \,\, (t,\theta) \in [\lambda,t_0] \times {\mathbb S}_{n-1}. $$

\end{Lemma}

\bigskip

Let $ \xi_k $ be the following real number,

$$ \xi_k=\sup \{ \lambda \leq \lambda_k+2, \bar w_k^{\lambda}(t,\theta)-\bar w_k(t,\theta) \leq 0, \,\, \forall \,\, (t,\theta)\in [\lambda,t_0]\times {\mathbb S}_{n-1} \}. $$

We have the same computations as for the previous papers, see [4,6,8,9,10]. We have the increment of functions and operators.

$$ \tilde L(t,\theta )[\tilde w_k^{\xi_k}(t, \theta)-\tilde w_k(t,\theta)]= [\tilde L(t,\theta )-\tilde L(t^{\xi_k},\theta)][\tilde w_k^{\xi_k}(t, \theta)]+ \tilde L(t^{\xi_k},\theta )[\tilde w_k^{\xi_k}(t, \theta)]- \tilde L(t,\theta )[\tilde w_k(t,\theta)]=$$

$$ =[\tilde L(t,\theta )-\tilde L(t^{\xi_k},\theta)][\tilde w_k^{\xi_k}(t, \theta)]+ \tilde V_k^{\xi_k} \left (\dfrac{1}{ \tilde b_1^{\xi_k}} \right )^{N-2} (\tilde w_k^{\xi_k})^{(n+2)/(n-2)}- \tilde V_k \left (\dfrac{1}{ \tilde b_1} \right )^{N-2} \tilde w_k^{(n+2)/(n-2)}$$

Thus,

$$\tilde L(t,\theta )[\tilde w_k^{\xi_k}(t, \theta)-\tilde w_k(t,\theta)] =  $$

$$ =[\tilde L(t,\theta )-\tilde L(t^{\xi_k},\theta)][\tilde w_k^{\xi_k}(t, \theta)]+ \tilde V_k^{\xi_k} \left (\dfrac{1}{ \tilde b_1^{\xi_k}} \right )^{N-2} (\tilde w_k^{\xi_k})^{(n+2)/(n-2)}- \tilde V_k \left (\dfrac{1}{ \tilde b_1} \right )^{N-2} \tilde w_k^{(n+2)/(n-2)}. $$

Thus,

$$\tilde L(t,\theta )[\bar w_k^{\xi_k}(t, \theta)-\bar w_k(t,\theta)] =  $$

$$ =[\tilde L(t,\theta )-\tilde L(t^{\xi_k},\theta)][\tilde w_k^{\xi_k}(t, \theta)]+ $$

$$ + \tilde V_k^{\xi_k} \left (\dfrac{1}{ \tilde b_1^{\xi_k}} \right )^{N-2} (\tilde w_k^{\xi_k})^{(n+2)/(n-2)}- \tilde V_k \left (\dfrac{1}{ \tilde b_1} \right )^{N-2} \tilde w_k^{(n+2)/(n-2)}+ $$

$$ + O(1) r_k^{4/(n-2)} e^{2t} (e^{(n-2)t/2}-e^{(n-2)t^{\xi_k}/2}). $$

We have:

$$ [\tilde L(t,\theta )-\tilde L(t^{\xi_k},\theta)][\tilde w_k^{\xi_k}(t, \theta)]=O(1) \tilde w_k^{\xi_k}r_k^{4/(n-2)} (e^{2t}-e^{2t^{\xi_k}}), $$

and,

$$ \tilde V_k^{\xi_k} \left (\dfrac{1}{ \tilde b_1^{\xi_k}} \right )^{N-2} (\tilde w_k^{\xi_k})^{(n+2)/(n-2)}- \tilde V_k \left (\dfrac{1}{ \tilde b_1} \right )^{N-2} \tilde w_k^{(n+2)/(n-2)}= $$

$$ = O(1) r_k^{4/(n-2)} \tilde w_k^{\xi_k} (e^{2t}-e^{2t^{\xi_k}}) +O(1)(\tilde w_k^{\xi_k})^{(n+2)/(n-2)} r_k^{2/(n-2)} (e^t-e^{t^{\xi_k}}) + $$

\be + \tilde V_k \left (\dfrac{1}{ \tilde b_1} \right )^{N-2} [(\tilde w_k^{\xi_k})^{(n+2)/(n-2)}- \tilde w_k^{(n+2)/(n-2)} ]. \ee

Thus,

$$\tilde L(t,\theta )[\bar w_k^{\xi_k}(t, \theta)-\bar w_k(t,\theta)] =  $$

$$ =\tilde V_k \left (\dfrac{1}{ \tilde b_1} \right )^{N-2} [(\tilde w_k^{\xi_k})^{(n+2)/(n-2)}- \tilde w_k^{(n+2)/(n-2)} ] + $$

$$ +O(1) r_k^{4/(n-2)} \tilde w_k^{\xi_k} (e^{2t}-e^{2t^{\xi_k}}) + $$

$$ + O(1)(\tilde w_k^{\xi_k})^{(n+2)/(n-2)} r_k^{2/(n-2)} (e^t-e^{t^{\xi_k}}) + $$

\be + O(1) r_k^{4/(n-2)} e^{2t} (e^{(n-2)t/2}-e^{(n-2)t^{\xi_k}/2}). \ee

We want to prove that by using the Hopf maximum principle, (like in [3,5,7,8]):

$$ \min_{\theta \in {\mathbb S}_{n-1}} \bar w_k(t_0,\theta) \leq \max_{\theta \in {\mathbb S}_{n-1}} \bar w_k(2\xi_k-t_0, \theta), $$

For this, we argue by contradiction and we assume that:

$$ \min_{\theta \in {\mathbb S}_{n-1}} \bar w_k(t_0,\theta) > \max_{\theta \in {\mathbb S}_{n-1}} \bar w_k(2\xi_k-t_0, \theta), $$

Thus, our assumption is:

$$ \bar w_k(2\xi_k-t_0, \theta)- w_k(t_0,\theta) <0, \forall \theta \in {\mathbb S}_{n-1}. $$

Now, we want to prove that:

$$ [\bar w_k^{\xi_k}(t, \theta)-\bar w_k(t,\theta)] \leq 0 \Rightarrow 
   \tilde L(t,\theta )[\bar w_k^{\xi_k}(t, \theta)-\bar w_k(t,\theta)] \leq 0, $$

For this:

\bigskip

1) The biggest term is the term of $ V $ (for $ n\geq 4$): $ \tilde w_k^{\xi_k} r_k^{2/(n-2)} (e^t-e^{t^{\xi_k}}), t_0\geq t\geq \xi_k $. 

\smallskip

Because we must compare:

$$ (\tilde w_k^{\xi_k})^{(n+2)/(n-2)} r_k^{2/(n-2)} (e^t-e^{t^{\xi_k}})\,\,{\rm and} \,\, (\tilde w_k^{\xi_k})^{(n+2)/(n-2)}-(\tilde w_k)^{(n+2)/(n-2)}, $$

and, we have used the mean value theorem for $ f(t)=t^{(n+2)/(n-2)} $, and $ \tilde w_k^{\xi_k}\leq \tilde w_k $ to have: 

$$(\tilde w_k^{\xi_k})^{(n+2)/(n-2)}-(\tilde w_k)^{(n+2)/(n-2)} \leq c (\tilde w_k^{\xi_k})^{4/(n-2)}(\tilde w_k^{\xi_k}-\tilde w_k ),  $$

and,

$$ (\tilde w_k^{\xi_k})^{4/(n-2)}(\tilde w_k^{\xi_k}-\tilde w_k )\leq c (\tilde w_k^{\xi_k})^{4/(n-2)}(e^{(n-2)t^{\xi_k}/2}-e^{(n-2)t/2}). $$

Now, we write:

$$ e^{t}=e^{(n-2)t/2} e^{(4-n)t/2} \leq e^{(4-n)\xi_k/2} e^{(n-2)t/2}, $$

we integrate between $ t $ and $ t^{\xi_k} $, we obtain:

$$ (e^{t}-e^{t^{\xi_k}})\leq c e^{(4-n)\xi_k/2} (e^{(n-2)t/2}-e^{(n-2)t^{\xi_k}/2}), $$

But,

$$ \tilde w_k^{\xi_k} \leq c e^{(n-2)(\xi_k-\lambda_k)/2}, $$

Thus the biggest term is:

$$ \tilde w_k^{\xi_k} r_k^{2/(n-2)} (e^t-e^{t^{\xi_k}}) \leq c r_k^{2/(n-2)}e^{(2\xi_k-(n-2)\lambda_k)/2} (e^{(n-2)t/2}-e^{(n-2)t^{\xi_k}/2}), $$

but $ \xi_k \leq \lambda_k+2 $, we obtain:

$$ \tilde w_k^{\xi_k} r_k^{2/(n-2)} (e^t-e^{t^{\xi_k}}) \leq c r_k^{2/(n-2)} e^{-(n-4)\lambda_k/2}(e^{(n-2)t/2}-e^{(n-2)t^{\xi_k}/2}), $$

Thus,

$$ n\geq 5, \,\, \tilde w_k^{\xi_k} r_k^{2/(n-2)} (e^t-e^{t^{\xi_k}}) \leq \frac{c}{[u_k(t_k)]^{\epsilon-(n-4)/(n-2)}} (e^{(n-2)t/2}-e^{(n-2)t^{\xi_k}/2}), $$

$$ n=4, \,\, \tilde w_k^{\xi_k} r_k^{2/(n-2)} (e^t-e^{t^{\xi_k}}) \leq \frac{c}{[u_k(y_k)]^{\epsilon}} (e^{(n-2)t/2}-e^{(n-2)t^{\xi_k}/2}), $$

These terms are controled by the term: $ -\frac{m}{2}(e^{(n-2)t/2}-e^{(n-2)t^{\xi_k}/2}) $.

\bigskip

2) Also, we have for $ n\geq 6 $:

\smallskip

Because we must compare:

$$ \tilde w_k^{\xi_k} r_k^{4/(n-2)} (e^{2t}-e^{2t^{\xi_k}})\,\,{\rm and} \,\, (\tilde w_k^{\xi_k})^{(n+2)/(n-2)}-(\tilde w_k)^{(n+2)/(n-2)}, $$

\smallskip

We must look to the term: 

$$ (\tilde w_k^{\xi_k})^{1-(4/(n-2))} r_k^{4/(n-2)} (e^{2t}-e^{2t^{\xi_k}}) =(\tilde w_k^{\xi_k})^{(n-6)/(n-2)} r_k^{4/(n-2)} (e^{2t}-e^{2t^{\xi_k}}), t_0\geq t\geq \xi_k, $$

We write:

$$ e^{2t}=e^{(n-2)t/2} e^{(6-n)t/2} \leq e^{(6-n)\xi_k/2} e^{(n-2)t/2}, $$

But,

$$ \tilde w_k^{\xi_k} \leq c e^{(n-2)(\xi_k-\lambda_k)/2}, $$

Thus,

$$ (\tilde w_k^{\xi_k})^{(n-6)/(n-2)} r_k^{4/(n-2)} (e^{2t}-e^{2t^{\xi_k}}) \leq c r_k^{4/(n-2)}e^{-(n-6)\lambda_k/2} (e^{(n-2)t/2}-e^{(n-2)t^{\xi_k}/2}), $$

we obtain:

$$ (\tilde w_k^{\xi_k})^{(n-6)/(n-2)} r_k^{4/(n-2)} (e^t-e^{t^{\xi_k}}) \leq c r_k^{4/(n-2)} e^{-(n-6)\lambda_k/2(}(e^{(n-2)t/2}-e^{(n-2)t^{\xi_k}/2}), $$

Thus,

$$ n\geq 6, \,\, (\tilde w_k^{\xi_k})^{(n-6)/(n-2)} r_k^{4/(n-2)} (e^{2t}-e^{2t^{\xi_k}}) \leq \frac{c}{[u_k(t_k)]^{\epsilon-(n-6)/(n-2)}} (e^{(n-2)t/2}-e^{(n-2)t^{\xi_k}/2}), $$

But, $ \epsilon >\frac{n-4}{n-2} $, these terms are controled by the term: $ -\frac{m}{2}(e^{(n-2)t/2}-e^{(n-2)t^{\xi_k}/2}) $.

\bigskip

3) For $ n=5 $: we have the terms: we use the binomial formula: we write:

$$ (\tilde w_k^{\xi_k})^{7/3}-\tilde w_k^{7/3}=((\tilde w_k^{\xi_k})^{1/3})^7-(\tilde w_k^{1/3})^7, $$

$$ x^7-y^7\equiv (x-y)(x^6+x^5y+x^4y^2+x^3y^3+x^2y^4+xy^5+y^6), x=(\tilde w_k^{\xi_k})^{1/3}, y=\tilde w_k^{1/3}, $$

but,

$$ \tilde w_k^{\xi_i}-\tilde w_k = (x^3-y^3) \equiv (x-y)(x^2+xy+y^2), $$

Thus,

$$ (\tilde w_k^{\xi_k})^{7/3}-\tilde w_k^{7/3}=(x^3-y^3) \times \frac{(x^6+x^5y+x^4y^2+x^3y^3+x^2y^4+xy^5+y^6)}{(x^2+xy+y^2)}, $$

Here, we have $ x\leq y $, thus:

$$ (x^2+xy+y^2) \leq c y^2, (x^6+x^5y+x^4y^2+x^3y^3+x^2y^4+xy^5+y^6) \geq c'x^2y^4, c, c'>0 $$

Thus, because $ y \geq \frac{m}{2} e^{3t/2} $ we obtain:

$$ (\tilde w_k^{\xi_k})^{7/3}-\tilde w_k^{7/3} \leq c(\tilde w_k^{\xi_k}-\tilde w_k)(\tilde w_k^{\xi_k})^{2/3} \tilde w_k^{2/3} \leq -c e^{t} (\tilde w_k^{\xi_k})^{2/3} (e^{3t/2}-e^{3t^{\xi_k}/2}), c>0 $$

For the case: $ A=(\tilde w_k^{\xi_k})^{1/3} r_k^{4/3} (e^{2t}-e^{2t^{\xi_k}}), t_0\geq t\geq \xi_k $

We have:

$$ |A|\leq e^t(\tilde w_k^{\xi_k})^{1/3} r_k^{4/3} (e^{t}-e^{t^{\xi_k}}) $$

The dominant term is:

$$ B=(\tilde w_k^{\xi_k})^{1/3} r_k^{4/3} (e^{t}-e^{t^{\xi_k}})$$

We have:

$$ e^t=e^{-t/2}e^{3t/2} \leq ce^{-\xi_k/2} e^{3t/2}, w_k^{\xi_k}\leq c e^{3(\xi_k-\lambda_k)/2} $$

Thus,

$$ |B|\leq c r_k^{4/3} e^{-\lambda_k/2} (e^{3t/2}-e^{3t^{\xi_k}/2})$$

$$ \lambda_k=-(2/3)(1-\epsilon) \log u_k(t_k), r_k=u_k(t_k)^{-\epsilon}, $$

$$ e^{-\lambda_k/2}=u_k(t_k)^{(1/3)(1-\epsilon)}, r_k^{4/3}=u_k(t_k)^{-4\epsilon/3}, r_k^{4/3}e^{-\lambda_k/2}= u_k(t_k)^{(1/3)(1-5\epsilon)}  $$

The condition is $ 1-5\epsilon <0 $, $ \epsilon > \frac{1}{5} $, but $ \epsilon >\frac{n-4}{n-2}=\frac{5-4}{5-2}=1/3 >1/5 $.

\bigskip

4) We have the same thing for the dimension 4.

\smallskip

5) When we use the auxiliary function $ \frac{m}{2} e^{(n-2)t/2} $, there is a term:

$$ r_k^{4/(n-2)} e^{2t} (e^{(n-2)t/2}-e^{(n-2)t^{\xi_k}/2}), $$ 

To correct this term, we consider a part of the term :

$$(\tilde w_k^{\xi_k})^{(n+2)/(n-2)}-(\tilde w_k)^{(n+2)/(n-2)} $$

We use the binomial formula as for the previous case of dimension 5. We have:

$$ x=(\tilde w_k^{\xi_k})^{1/(n-2)}, y= (\tilde w_k)^{1/(n-2)}, $$

$$ x^{n+2}-y^{n+2}=(x-y) (y^{n+1}+\ldots) , x^{n-2}-y^{n-2}=(x-y)(y^{n-3}+\ldots+x^{n-3}) $$

Thus,

$$ x^{n+2}-y^{n+2}=(x^{n-2}-y^{n-2})\frac{(y^{n+1}+\ldots)}{(y^{n-3}+\ldots)}$$

Because $ x\leq y $ and,

$$ (y^{n+1}+\ldots) \geq y^{n+1}, (x^{n-3}+\ldots+y^{n-3}) \leq cy^{n-3}, c >0,$$

We obtain:

$$ x^{n+2}-y^{n+2} \leq c(x^{n-2}-y^{n-2}) y^4, c>0,$$

Thus,

$$(\tilde w_k^{\xi_k})^{(n+2)/(n-2)}-(\tilde w_k)^{(n+2)/(n-2)} \leq c( \tilde w_k^{\xi_k}-\tilde w_k ) \tilde w_k^{4/(n-2)}, c >0 $$

Because, $ \tilde w_k\geq \frac{m}{4} e^{(n-2)t/2} $, we obtain:

$$(\tilde w_k^{\xi_k})^{(n+2)/(n-2)}-(\tilde w_k)^{(n+2)/(n-2)} \leq -c e^{2t} (e^{(n-2)t/2}-e^{(n-2)t^{\xi_k}/2}), c >0, $$

Thus the term: $ r_k^{4/(n-2)} e^{2t} (e^{(n-2)t/2}-e^{(n-2)t^{\xi_k}/2}) $ is controled by the term, $ -c e^{2t} (e^{(n-2)t/2}-e^{(n-2)t^{\xi_k}/2}), c >0 $.

\bigskip

We obtain the same proof in the previous papers, the dimensions $ 4,6 $, see [6,9,10].

\bigskip

If we use the Hopf maximum principle, we obtain (like in [4,6,8,9]):

$$ \min_{\theta \in {\mathbb S}_{n-1}} \bar w_k(t_0,\theta) \leq \max_{\theta \in {\mathbb S}_{n-1}} \bar w_k(2\xi_k-t_0, \theta), $$

thus for $ k $ large:

$$ +\infty \leftarrow (u_k(\cdot))^{1-\epsilon}=v_k(0) \leq c, $$

it is a contradiction.

\bigskip

Finaly, for each $ m >0 $ there is a finite $ v_k $ such that $ \inf v_k \geq m >0, k_1,\ldots, k_m \in {\mathbb N} $. 

\bigskip

Here also, we have the existence of $ c(m) >0 $ such that $ \inf v_k \geq m >0 \Rightarrow (u_k(\cdot))^{1-\epsilon} = v_k(0) \leq c $. We prove this by contradiction, suppose that for fixed $ m >0 $, for all $ c >0 $ there is $ i_c\in {\mathbb N} $ with $ \inf v_{i_c} \geq m >0 $ and $ v_{i_c}(0) \geq c $, if we take $ c \to +\infty $, because the number of indices is bounded and we have a sequence of integers, this sequence converge and in fact is constant because we consider integers. Thus there is an index $ k $ such that $ v_k(0) \geq c \to +\infty $ and $ \inf v_k \geq m >0 $ , and thus $ v_k $ is singular at $ 0 $, but this is impossible because $ v_k $ is regular.

\smallskip

We obtain:

\smallskip

There is a non-increasing positive function $ m \to c(m)>0 $, such that $ \inf v_k \geq m >0 \Rightarrow  (u_k(\cdot))^{1-\epsilon} \leq c(m) $. then we apply this with $ m = r_k \inf_M u_k $, we obtain the inequality for all terms of the sequence $ (u_k)$.

\bigskip

\bigskip

{\bf Remark: Compactness Criterion and Compactness result:}

\bigskip

We have the following compactness criterion and compactness result: the sequence $ (u_k) $ is compact around the point $ x_0 $, if and only if: it satisfies locally the usual Harnack inequality around the point $ x_0 $.

\bigskip

The local usual Harnack inequality around the point $ x_0 $: 

$$ \forall \,\, (x_k), (y_k), \,\, x_k \to x_0, y_k \to x_0,\,\, \exists \, m >0, u_k(x_k)\geq m \cdot u_k(y_k), \forall k. \,\,\, (**) $$

We have:

$$ (**) \Leftrightarrow {\rm \, compactness \,\, around\,\,}  x_0 .$$

Indeed, we can apply (up to a subsequence) this definition to the points $ x_k $ of the infimum $ \inf_{B(t_k,2r_k^{2/(n-2)})} u_k = u_k(x_k) $ and $  y_k=t_k $, in the previous proof of the Theorem 1.1, then, the point  2) of the Theorem 1.1. is not possible. Thus, we have the point 1) of the Theorem 1.1, and then, we have the compactness result around $ x_0 $.

\bigskip

Here, we choose the radius of the ball, small and tending to 0, to have not, the uniform boundedness of the infimum, as showed in the paper of Brezis-Merle. If we take a fixed radius for the infimum, this infimum become bounded uniformly, by the argument of the first eigenvalue of Brezis-Merle, see [12], and there is nothing to prove.

\bigskip

This is a reason why we take this definition of local usual Harnack inequality: around the points and in particular $ x_0 $.

Here we talk about a compactness criterion and also about a compactness result, because, the usual example of functions: $ x\to (\epsilon/(\epsilon^2+|x|^2))^{(n-2)/2}, \epsilon \to 0 $, do not satisfy the local usual Harnack inequality and $ (**) $ at 0. Thus, to have the compactness result in all generality, we must suppose $ (**) $.

\end{document}